\documentclass[journal]{IEEEtran}
\usepackage{}

\usepackage{graphicx}
\usepackage{amsmath, amssymb}

\newtheorem{theorem}{Theorem}
\newtheorem{lemma}{Lemma}
\newtheorem{remark}{Remark}
\newtheorem{proposition}{Proposition}

\newtheorem{example}{Example}
\usepackage{CJK}
\newcommand{\be}{\begin{equation}}
\newcommand{\ee}{\end{equation}}
\ifCLASSINFOpdf
\else
\fi
\begin{document}
%
\title{Learning Based Control Policy and Regret Analysis for Online Quadratic Optimization with Asymmetric Information Structure}
%
%
%

\author{Cheng~Tan~\IEEEmembership{Member,~IEEE} and Wing~Shing~Wong~\IEEEmembership{Fellow,~IEEE}
\thanks{This work was supported in part by a grand from the Research Grants Council of the Hong Kong Special Administrative Region under Project GRF. 14630915,
the Taishan Scholar Project of Shandong Province of China ts201712040,
the National Natural Science Foundation of China under Grants 61803224,
the Natural Science Foundation of Shandong Province ZR. 201702170323.
}

\thanks{C. Tan and W. S. Wong are with the Department of Information Engineering, The Chinese University of Hong Kong, Shatin, N. T., Hong Kong (\small  e-mail: tancheng1987love@163.com; wswong@ie.cuhk.edu.hk).
C. Tan is also with the College of Engineering, QuFu Normal University, Rizhao, Shandong 276800, China.}
}
\maketitle

\begin{abstract}
In this paper, we propose a learning approach to analyze dynamic systems with asymmetric information structure. Instead of adopting a game theoretic setting, we investigate an online quadratic optimization problem driven by system noises with unknown statistics. Due to information asymmetry, it is infeasible to use classic Kalman filter nor optimal control strategies for such systems. It is necessary and beneficial to develop a robust approach that learns the probability statistics as time goes forward. Motivated by online convex optimization (OCO) theory, we introduce the notion of regret, which is defined as the cumulative performance loss difference between the optimal offline known statistics cost and the optimal online unknown statistics cost. By utilizing dynamic programming and linear minimum mean square biased estimate (LMMSUE), we propose a new type of online state feedback control policies and characterize the behavior of regret in finite time regime. The regret is shown to be sub-linear and bounded by $O(\ln T)$. Moreover, we address an online optimization problem with output feedback control policies.
\end{abstract}


%
\IEEEpeerreviewmaketitle

\section{Introduction}

Many previously reported works on dynamic systems assume the classic information structure that postulates all agents have equal access to available system information.
Such a symmetric information structure is encountered in a host of application scenarios such as pursuit-evasion games \cite{pe01}-\cite{pe02}, networked control systems \cite{wong1999}-\cite{tan2017} and seller-buyer supply chain models \cite{sb01}-\cite{sb02}.
In differential game settings, it is common to assume that the opposing parties have peering information in regard to location, velocity, player utility functions and control policies.
For example in \cite{sb01}, the seller and the buyer achieve the pricing and batch-size equilibrium by solving a cooperative Stackelberg game.
While such a symmetric information assumption is satisfied in many applications, from a general application perspective it is of interest to investigate systems with an asymmetric information structure.
Moreover, early pioneering work in \cite{as01}-\cite{as02} has pointed out the important role played by the information structure on decision and control strategy, and thus offering theoretical motivation to study systems with a non-classic information structure.
There are a number of works analyzing models with asymmetric information in dynamic games \cite{dy01}-\cite{dy03},  pursuit-evasion \cite{pe03}, and economic theory  \cite{sb03}-\cite{sb05}.

In this paper, we aim to analyze two-player systems in which a single agent with rich input information, the {predator}, is pitted against the other agent with limited input information, the {prey}.
The motivation of the model comes from application scenarios that include pursuit-evasion and product pricing.
Due to its asymmetric nature, we formulate the problem as a quadratic optimization from the perspective of the {predator} instead of a game theoretic setting.
Below, we use two simple examples to illustrate the types of online quadratic optimization we focus on in this paper.

The first example is based on the pursuit-evasion model in \cite{taniet2017}
and the Mission 7 challenge of the International Aerial Robotics Competition (IARC) in \cite{iarc},
consisting of a single predator and a single prey.  The predator has access to location information of both players and based on that selects a predation mode (for example whether bait or camouflage is used) and a predation policy at each decision instant.
The prey is aware of the selected predation mode but otherwise has no access to location information of the predator.
Hence, it adopts a simple randomized evading policy for each predation mode.
A simple illustration of the dynamic game is depicted in Fig. 1, where the blue lines represent trajectories corresponding to the 1st predation mode and the red lines represent trajectories of the 2nd predation mode.
\vspace{-5mm}
\begin{figure}[thpb]
      \centering
      \includegraphics[scale=0.55]{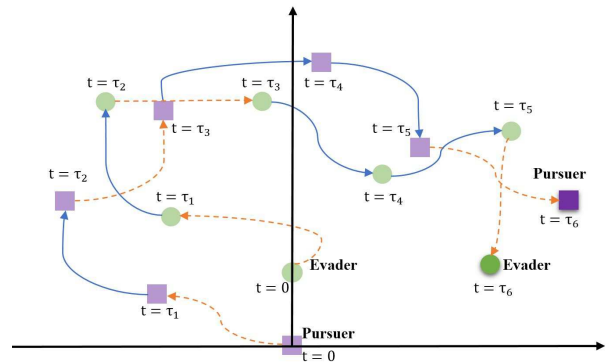}
      \caption{Movement trajectories of pursuer and evader \cite{taniet2017}}
\end{figure}

To be specific, the dynamic of the predator and the prey is described as
\begin{align}
x_p(t+1)=&~x_p(t)+u(i(t),t), \label{sys00a}\\
x_e(t+1)=&~x_e(t)+v_{i(t)}(t), \label{sys00b}
\end{align}
where $x_p(t)$, $x_e(t)$ are the respective positions of the predator and prey.
Denote $i(t)\in\mathbb{K} \triangleq\{1,2,\cdots,K\}$ to be the predation mode and $u(i(t),t)$ the predation policy, which are to be determined.
Assume the evading policies are defined by random variables, $v_{i(t)}(t)$'s, which are independent of each other.
Moreover, for each mode $k\in\mathbb{K}$, $v_{k}(t)$ takes value in an admissible bounded set $\{v_1,\cdots,v_M\}$ with $prob(v_k(t)=v_i)=p_{k,i}$.
The objective of the predator is to minimize both control cost and distance, which is captured by the following index function
\begin{equation}\label{cost00}
W_T
=\sum_{t=0}^T \beta^t{\bf E}
\left[\|x_e(t)-x_p(t)\|^2+\left\|u(i(t),t)\right\|^2\right],
\end{equation}
with a given $0<\beta\leq 1$.
We emphasize that the evading policy distributions are a priori unknown to the predator.

Our second example is related to product pricing \cite{price}.
Consider a product pricing that is determined by a single producer, which has absolute control over the pricing and the producing rate.
The market demand rate of the product, $d(t)$, satisfies the model
\begin{equation}
d(t+1)=\max\{\zeta(t)-\alpha p(t),0 \},
\end{equation}
where $\alpha>0$ and $p(t)$ is the pricing set by the producer.
$\zeta(t)$ represents the utility value of the product and is assumed to satisfy
\begin{equation}
\zeta(t)=b+e(t),
\end{equation}
where $e(t)$'s are independent and identically distributed (i.i.d.) random variables with zero mean and variance $v_e$.
If $\zeta(t)$ is assumed to be positive and bounded away from zero and $\alpha>0$  is relatively
small, we simplify the demand model as
\begin{equation}
d(t+1)=\zeta(t)-\alpha p(t),~d(0)=d_0.\label{price01}
\end{equation}
The production process is modelled by
\begin{equation}
z(t+1)=z(t)+u(t),\label{price02}
\end{equation}
where $z(t)$ is the production rate and $u(t)$ is the rate control.
Moreover, $z(0)$ is a given constant, $z_0$, known to the producer.
For a given optimization horizon of $T$ periods, we define the following objective function
\begin{align}
J_T = \sum_{t=0}^T {\bf{E}}& \big[ c_1  (z(t)-d(t))^2 + c_2 u^2 (t)  \nonumber\\
&  - c_3 (p(t)-C)d(t+1) \big] , \label{opti01}
\end{align}
where $c_i>0,~i=1,2,3$, are known positive weighting constants and $C>0$ is the product cost.
The first component in $J_T$ measures how the production process tracks the demands;
the second term is a measure of the production rate changes;
the last component  represents the total profit assuming the demands are met.
The objective of the producer is to minimize $J_T$ via the control variables $(u(t), p(t))$,
which are assumed to be measurable with respect to (w.r.t.) the $\sigma$-algebra generated
by the set $\{z(s),d(s),s=0,1,\cdots,t\}$.

Suppose $b>0$ is known to the producer.
If we set
\begin{equation}
v(t)=-p(t)+\frac{1}{2}(C+\frac{b}{\alpha}),~y(t)=\frac{d(t)}{\alpha},
\end{equation}
the objective function (\ref{opti01}) can be reformulated as
\begin{align}
J_T =&~\bar{J}_T +\sum_{t=0}^T {\bf{E}}\left[ bC-\frac{1}{4}(C+\frac{\zeta(t)}{\alpha})^2  \right], \label{opti01a}
\end{align}
where
\begin{align}
&\bar{J}_T=\sum_{t=0}^T  {\bf{E}} \big[ c_1 (z(t)-\alpha y(t))^2
+ c_2 u^2 (t)+ c_3 \alpha v^2 (t) \big], \label{opti01b} \\
&~~~~~~~~~y(t+1)=v(t)+w(t),  \label{dy01}  \\
&~~~~~~~~~z(t+1)=z(t)+u(t),\label{dy02} \\
&~~~~~~~~~~~~~w(t)=\frac{e(t)}{\alpha}
-\frac{1}{2} (C-\frac{b}{\alpha}). \label{dy02a}
\end{align}
In this case,  the control variables $(u(t), v(t))$ are measurable w.r.t. the $\sigma$-algebra generated by $\{z(s),~y(s),~s=0,1,\cdots,t\}$.
The difference between the two objective functions $J_T$ and $\bar{J}_T$ is independent of the control policy.
Hence, the original problem can be reduced to minimize the quadratic function $\bar{J}_T$ subject to (\ref{dy01})-(\ref{dy02a}).

In the two simple examples above, due to the asymmetric information structure, the probability statistics of $v_{i(t)}(t)$ in  (\ref{sys00b}) and $e(t)$ in (\ref{dy02a}) are a priori unknown to the predator and the producer respectively.
Therefore, due to information asymmetry, it is infeasible to use the well known {dynamic programming approach} \cite{bert1995} nor {the maximum principle} for such systems \cite{ery1975}-\cite{wuzhen}.
It is necessary as well as beneficial to develop a robust approach that learns as time goes forward.

The framework of online convex optimization (OCO) was first defined in the machine learning literature \cite{zink2003}-\cite{sha2018}, which is closely tied to statistical learning theory and convex optimization.
In OCO theory, an online player iteratively makes decisions, whose ultimate goal is to minimize the cumulative cost in a long run which translates to making fewer prediction mistakes in the classification case.
The popular performance metric for online algorithms is regret.
In principle, the regret analysis aims to study how far an online algorithm deviates from the optimum \cite{bubeck2012}.
To be specific, the {regret} is defined as the cumulative performance loss difference between the online cumulative unknown statistics cost $J_T(u)$ and the optimal offline known statistics cost $J^*_T$.
An important property is that the regret of an online algorithm grows at a sub-linear rate, which means the time average of the index function converges to the optimal value as $T$ approaches infinity, i.e., $\lim_{T\rightarrow\infty}(J_T(u)-J^*_T)/T=0$.
In OCO framework, various online algorithms have been proposed to attain a regret of $O(\sqrt{T})$, such as the online gradient decent method \cite{zink2003}-\cite{KV2005}, the stochastic gradient decent method \cite{hazan2014} and the online Newton step method \cite{boyd2010}.
In \cite{hazan2007}, when the cost function is strictly convex, the regret can be improved to $O(\ln T)$.

In our previous work \cite{taniet2017}, motivated by recent OCO methodology \cite{hazan2016}-\cite{pater2017}, we reformulate the first pursuit-evasion model above as a Multi-Armed Bandit problem.
Our objective is to find the balance between staying with the predation mode with lowest cost and exploring new options with might lower cost in the future.
The proposed Gittins Index based control policy can be computed based on a forward induction.
Although the proposed algorithm outperforms a random decision policy, its regret is proved to be linear.
How to improve the Gittins Index based control policy to ensure a sub-linear regret is challenging and remains an open question.

In this paper, we focus on two-player systems in which the players have asymmetric ability to information as motivated by the above examples.
Instead of adopting a game theoretic setting, we investigate the quadratic optimization  problem on the predator side.
In general, our learning based research methodology contains three powerful techniques, namely, {dynamic programming}, {linear minimum mean square unbiased estimate (LMMSUE)} and {regret analysis}.

We formulate the problem as an online quadratic optimization driven by system noises with unknown statistics.
For the state feedback case, if the mean and variance of the system noises are known a priori, the optimal offline control policy is derived based on the dynamic programming approach.
The optimal state feedback gains, independent of the unknown statistics, are uniquely determined by solving a standard Riccati equation.
However, in the current model, since the probability statistics of the system noises are unknown, it is infeasible to apply the optimal offline control strategies.
To address this, we introduce a robust approach that learns the probability statistics of the system noises with the LMMSUE.
Based on that we propose a learning based optimal control policy.
Moreover, under some basic assumptions, the regret of the proposed online control policy grows at a sub-linear rate, which is shown to be bounded by $O(\ln T)$.
Simulation results show the performance of the developed control policy.
On the other hand, we try to address the online quadratic optimization problem with output feedback control.
Due to information asymmetry, the classic Kalman filter cannot be applied directly.
With the LMMSUE, we propose a heuristic online control policy.
The regret between the online known statistics cost and the proposed heuristic offline unknown statistics cost is sub-linear, that shown to be bounded by $O(\ln T)$.

The following is an outline of this paper.
In Section II we investigate an online quadratic optimization problem with state feedback control. A LMMSUE-based online control policy is developed whose regret is shown to be bounded by $O(\ln T)$.
In Section III we address an output feedback case.
In Section IV, two simple examples are presented to illustrate the effectiveness of the developed control policies.
A conclusion is presented in Section V along with some relevant remarks.

{\em Notation}:
Let $\mathbb{R}^n$ denote the $n$-dimensional real Euclidean space and $\mathbb{R}^{m\times n}$ be the space formed by all $m\times n$ real matrices with the usual 2-norm $\| \cdot \|$.
The superscript $'$ represents matrix transpose.
${\bf Tr}(A)$ represents the trace of a square matrix $A$ and $diag\{a_1~a_2~\cdots~a_n\}$ denotes a diagonal matrix.
$A\geq 0~(>0)$ represents that $A$ is a positive semi-definite (positive definite) matrix and $A\geq B\ (>B)$ means that $A-B\geq 0~(>0)$.
$\{w(t),t=0,1,\cdots\}$ denotes a sequence of real random variables defined on the
complete filtered probability space $(\Omega,
\mathcal{F},\mathcal{F}_t)$ with $\mathcal{F}_0=\{\emptyset,\Omega\}$
and $\mathcal{F}_t=\sigma\{w(s)|s=0,1,2,\cdots,t\}$.
Moreover, $prob(A)$ denotes the probability if the event $A$ occurs and ${\bf E} [w(t)]$ the expectation of the random variable $w(t)$.

\section{State Feedback Control with Learning}

\subsection{Problem Formulation}
Consider the following discrete time dynamic system
\begin{equation}\label{sys01}
x(t+1)=Ax(t)+Bu(t)+w(t),
\end{equation}
where $x(t) \in\mathbb{R}^n$ is the state and $u(t) \in\mathbb{R}^m$ is the input control.
$A,~B$ are the known system parameters with the compatible dimensions and $x(0)=x_0 \in\mathbb{R}^n$ is the given initial state.
We assume that $w(t)$'s, are bounded and i.i.d. stochastic process with
\begin{align}
&prob(w(t)=w_i)=p_i,~i=1,2,\cdots,M, \label{wt01}\\
&\max_{i}\|w_i\|\leq w_b<\infty. \label{wt02}
\end{align}
Define $\mathbf{p}_w=[p_1~p_2~\cdots~p_M]'$, $\mathbf{P}_w=diag\{p_1~p_2~\cdots~p_M\}$, and $\mathbf{W}=[w_1~w_2~\cdots~w_M]$.
It follows that
\begin{align}
\mu_w=&~{\bf E}[w(t)]=\sum_{i=1}^M p_iw_i=\mathbf{W}\mathbf{p}_w, \\
Q_w=&~{\bf E}[w(t)w(t)']=\sum_{i=1}^M p_iw_i w'_i=\mathbf{W}\mathbf{P}_w\mathbf{W}'.
\end{align}
Moreover, the covariance of $w(t)$ is
\begin{equation}
C_w= {\bf E} \big[(w(t)-\mu_w)(w(t)-\mu_w)'\big]=Q_w-\mu_w\mu_w'.
\end{equation}
Therefore, the probability statistic of $w(t)$ depends on $\mathbf{p}_w$.
We emphasize that $\mathbf{p}_w$ is a priori unknown to the decision maker, which leads to the asymmetric information structure.

Without loss of generality, the index function is defined as the general quadratic form
\begin{align}
J_T(u(t))= &\sum_{t=0}^{T} {\bf E} \big[ x'(t)Q(t)x(t)+u'(t)R(t)u(t) \big] \nonumber \\
&+{\bf E} \big[ x'(T+1)P_{T+1}x(T+1)\big], \label{cost01}
\end{align}
where $Q(t)\geq 0$, $R(t)>0$, and $P_{T+1}\geq 0$.
The goal of the decision maker is to minimize the index function (\ref{cost01}) by an online algorithm.

Suppose the probability $\mathbf{p}_w$ is known a priori. The finite horizon quadratic optimization problem (\ref{cost01}) subject to (\ref{sys01}) is fairly standard, which can be solved by utilizing the classic dynamic programming approach; see {Theorem \ref{thm01}} hereinafter.
Unfortunately, in the current model, $\mathbf{p}_w$ is unknown and the optimal known statistics control strategies cannot be applied directly for asymmetric information case.
How to address this unknown statistics problem?

Motivated by the OCO theory, we introduce the regret function as follows
\begin{align}
Reg_T(u(t))=J_T(u(t))-J_T^*. \label{regret00}
\end{align}
The regret measures the cumulative performance loss between the optimal offline case with known statistics cost $J_T^*$ and the online case with unknown statistics cost $J_T(u(t))$.
We say an online control policy performs well if its regret is sub-linear, i.e., $o(T)$, which implies the instantaneous online performance can converge asymptotically to that of the offline performance.
Our goal in this paper is to develop a robust approach to estimate the probability $\mathbf{p}_w$ based on the observed state trajectory and then propose a learning based control policy to reach a sub-linear regret.

\begin{remark}\label{remark01}
In the predator-prey model with single predation mode, i.e., $i(t)\equiv 1$, if we set $x(t)=x_p(t)-x_e(t)$, the first example can be equivalently reformulated as
\begin{align*}
&minimize~~J_T= \sum_{t=0}^{T} \beta^t {\bf E} \big[ x'(t)x(t)+u'(t)u(t) \big], \\
&subject~to~~x(t+1)=x(t)+u(t)-v(t),
\end{align*}
where $v(t)$ takes value in an admissible bounded set $\{v_1,\cdots,v_M\}$ with $prob(v_k(t)=v_i)=p_{i}$.
In the second product pricing example, if we set $X(t)=[y(t)~z(t)]'$, $U(t)=[v(t)~u(t)]'$ and $W(t)=[w(t)~0]'$, the original problem can be reformulated as
\begin{align*}
&minimize~~\bar{J}_T= \sum_{t=0}^{T} {\bf E} \big[ X'(t)QX(t)+U'(t)RU(t) \big], \\
&subject~to~~Z(t+1)=AX(t)+BU(t)+W(t),
\end{align*}
where
\[ \begin{matrix} A=\begin{bmatrix}
 0 & 0 \\ 0& 1
\end{bmatrix},
~B=\begin{bmatrix}
 1 & 0 \\ 0& 1
\end{bmatrix},
\end{matrix}\]
\[ \begin{matrix}
Q=c_1\begin{bmatrix}
 \alpha^2 & -\alpha \\ -\alpha & 1
\end{bmatrix}\geq 0,
~R=\begin{bmatrix}
 c_3 & 0 \\ 0 & c_2
\end{bmatrix}>0.
\end{matrix}\]
Therefore, both examples are the special cases of the online quadratic optimization problem with asymmetric information structure.
\end{remark}

\subsection{Preparatory Results}

To begin with, we derive the optimal control strategy $u^*(t)$ and the known statistics optimum $J_T^*$ based on perfect information of the probability $\mathbf{p}_w$.

\begin{theorem}
\label{thm01}
Suppose the probability $\mathbf{p}_w$ is known a priori.
The optimal offline control policy of the quadratic optimization problem (\ref{cost01}) is
\begin{align}
u^*(t)=&-\Upsilon_T(t)^{-1} \big( B'P_T(t+1)A x(t)+B'P_T(t+1)\mu_w \nonumber\\
&+B'L_T(t+1)\mu_w \big), \label{input}
\end{align}
while the optimal offline known statistics index value of (\ref{cost01}) is
\begin{equation}\label{cost02}
J_T^*=x_0'P_T(0)x_0+2x_0'L_T(0)\mu_w+\sum_{t=0}^TH_T(t),
\end{equation}
where $\Upsilon_T(t)$, $P_T(t)$, $L_T(t)$, and $H_T(t)$ satisfy the following iterative equations
\begin{align}
\Upsilon_T(t)=&~R(t)+B'P_T(t+1)B,  \label{ric00}\\
P_T(t)=&~ A'P_T(t+1)A+Q(t)-A'P_T(t+1)B\Upsilon_T(t)^{-1} \nonumber \\
&\times  B'P_T(t+1)A,\label{ric01} \\
L_T(t)=&~(A'-A'P_T(t+1)B \Upsilon_T(t)^{-1} B')(P_T(t+1)  \nonumber  \\
&+L_T(t+1)) \label{ric02}, \\
H_T(t)=&-\mu_w'(P_T(t+1)+L_T(t+1))'B \Upsilon_T(t)^{-1}B' \nonumber \\
&\times(P_T(t+1)+L_T(t+1))\mu_w \nonumber  \\
&+2\mu_w' L_T(t+1)\mu_w+{\bf Tr}\left(P_T(t+1)Q_w \right), \label{ric03}
\end{align}
with the terminal condition $P_T(T+1)=P_{T+1}$ and $L_T(T+1)=0$.
\end{theorem}

\textbf{Proof.} See Appendix~\ref{app01}. \hfill $\Box$

Since the probability $\mathbf{p}_w$ in the optimal offline control strategy is unknown a priori, the exact values of $\mu_w$ and $Q_w$ are unavailable.
Moreover, the optimum $J^*_T$ in Theorem \ref{thm01} is unavailable and can only be viewed as the optimal known statistics (offline) cost.

Note that $P_T(t)$, $L_T(t)$ are independent of $\mathbf{p}_w$ and thus can be computed offline.
Therefore,  $P_T(t)$, $L_T(t)$  are available to the decision maker at the initial time.
If we set $K_{P}(t)=-\Upsilon_T(t)^{-1} B'P_T(t+1)A$, the iterative Riccati equation (\ref{ric01}) is reduced to
\begin{align}
P_T(t)=&~ (A+BK_{P}(t))'P_T(t+1)(A+BK_{P}(t))\nonumber \\
&+Q(t)+K_{P}(t)'R(t)K_{P}(t). \label{ric01a}
\end{align}
Since $Q(t)\geq 0$ and $R(t)>0$, it follows from (\ref{ric01a}) that for any terminal condition $P_T(T+1)=P_{T+1}\geq0$, $P_T(t)\geq 0$ is unique and bounded.
Denote
\begin{equation}
\label{omega01}
\Omega_T(t)=A'-A'P_T(t+1)B \Upsilon_T(t)^{-1} B'.
\end{equation}
With the terminal condition $L_T(T+1)=0$, the adjoint equation (\ref{ric02}) can be rewritten as
\begin{align}
L_T(t)=&~\sum_{i=t+1}^{T+1}(\prod_{j=t}^{i-1}\Omega_T(j) ) P_T(i),\label{ric02b}
\end{align}
which indicates that the adjoint parameter $L_T(t)$ is uniquely determined by $P_T(s),~s=t+1,\cdots, T+1$ and is thus bounded.

Next, we evaluate the cost value in (\ref{cost01}) associated with any available control policy.

\begin{proposition}\label{prop01}
For any admissible control policy $u(t)$, the cost of the index function (\ref{cost01}) is
\begin{align}
&J_T(u(t)) \nonumber\\
=&~J_T^*+\sum_{t=0}^T {\bf E}\left( (u(t)-u^*(t))'\Upsilon_T(t) (u(t)-u^*(t)) \right), \label{eq011}
\end{align}
where $u^*(t)$ and $J_T^*$ are given in (\ref{input})-(\ref{cost02}).
\end{proposition}

\textbf{Proof.} See Appendix~\ref{app01a}. \hfill $\Box$

By Theorem \ref{thm01} and Proposition \ref{prop01}, we have
\begin{align*}
J_T^*-J_T(u(t))=\hspace{-0.8mm}\sum_{t=0}^T {\bf E}\big[ (u(t)-u^*(t))'\Upsilon_T(t) (u(t)-u^*(t))\big].
\end{align*}
For any admissible control policy $u(t)$, it follows from (\ref{regret00}) that the regret in this model can be rewritten as
\begin{align}
&Reg_T(u(t)) \nonumber \\
=&~\sum_{t=0}^T {\bf E}\big[ (u(t)-u^*(t))'\Upsilon_T(t) (u(t)-u^*(t))\big]. \label{error01}
\end{align}
For each time $t=1,2,\cdots,T$, and any admissible control policy $u(t)$, we define the one-step regret
\begin{equation}\label{error02}
reg_T(t,u(t))={\bf E} \big[ (u(t)-u^*(t))'\Upsilon_T(t) (u(t)-u^*(t)) \big].
\end{equation}
It follows that $Reg_T(u(t))=\sum_{t=0}^T reg_T(t,u(t))$.
The original optimization problem (\ref{cost01}) can be reduced to a minimization of (\ref{error01}) with some admissible online control policy.

\subsection{Learning Based Control Policy and Regret Analysis}

First, we focus on a simple but powerful learning tool of estimating an unknown parameter in statistical learning theory, i.e., LMMSUE.

Denote $\widehat{\mathbf{p}}_w(t)=[\hat{p}_{1}(t)~\hat{p}_{2}(t)~\cdots~\hat{p}_{M}(t) ]'$ to be the linear unbiased estimate of the probability $\mathbf{p}_w$.
With the initial estimate $\widehat{\mathbf{p}}_w(0)=[0~0~\cdots~0]'$,
it follows from \cite{taniet2017} that $\widehat{\mathbf{p}}_w(t)$ satisfies
\begin{equation}
\label{est}
\widehat{\mathbf{p}}_w(t)=\sum_{i=0}^{t-1}c_i(t)\xi(i),~t=1,2,\cdots,T,
\end{equation}
where $\sum_{i=0}^{t-1}c_i(t)=1$ and $\xi(i)$ is an i.i.d. stochastic process with
$prob(\xi(i)=\xi_j)=p_{j},~j=1,2,\cdots,M$,
and
$$\xi_1=[1~0~\cdots~0]', \cdots, \xi_M=[0~0~\cdots~1]'.$$
Actually, $\xi(t)$  defines the random observation that $w(t)$ takes the value of $w_i$ with the probability $p_i,~i=1,2,\cdots,M$.
In this case, we obtain that
\begin{align}
{\bf E}[\xi(i)]=\mathbf{p}_w,~{\bf E}[\xi(i)\xi(i)']=\mathbf{P}_w.
\end{align}
With linear unbiased estimate, we define the following admissible control policy set by
\begin{align}
\mathcal {U}_{ad}\triangleq \{u(t)=\hspace{-0.8mm}-\Upsilon_T(t)^{-1} B'P_T(t+1)Ax(t)+l_w(t)\},\label{controlset}
\end{align}
where
\begin{align*}
l_w(t)= -\Upsilon_T(t)^{-1} B' (P_T(t+1)+B' L_T(t+1)) \mathbf{W} \widehat{\mathbf{p}}_w(t).
\end{align*}
To begin with, we propose the LMMSUE $\widehat{\mathbf{p}}_{\min}(t)$ to minimize ${\bf E}\| \widehat{\mathbf{p}}_w(t)-\mathbf{p}_w\|^2$.

\begin{lemma}\label{lemma01}{\cite{taniet2017}}
The linear minimum mean square unbiased estimate of $\mathbf{p}_w$ is
\begin{equation}\label{est001}
\widehat{\mathbf{p}}_{\min}(t)=\frac{1}{t}\sum_{i=0}^{t-1} \xi(i).
\end{equation}
\end{lemma}

Note that the LMMSUE is the sample mean of the random observation $\xi(t)$.
It follows from (\ref{est001}) that
\begin{equation}
\widehat{\mathbf{p}}_{\min}(t+1)=\frac{1}{t+1}\left(t\widehat{\mathbf{p}}_{\min}(t)+\xi(t)\right).
\end{equation}
By the Kolmogorov Strong Law of Large Numbers \cite{ju1985}, we obtain
\begin{equation}
\lim_{t\rightarrow\infty}\widehat{\mathbf{p}}_{\min}(t)=\mathbf{p}_w,~a.s.
\end{equation}
where `a.s.' refers to `almost surely'.

\begin{remark}\label{remark02}
In principle, at each time $t=1,2,\cdots,T$, since $x(t)$, $x(t-1)$ and $u(t-1)$ are known to the decision maker, it is feasible to reach $w(t-1)$  with
\begin{equation}
w(t-1)=x(t)-Ax(t-1)-Bu(t-1).
\end{equation}
Observe that $w(t-1)=w_{h(t-1)}$,~$h(t-1)\in\mathbb{M} \triangleq\{1,2,\cdots,M\}$.
We update the LMMSUE
 $\widehat{\mathbf{p}}_{\min}(t)=(\hat{p}_{1}(t)~\hat{p}_{2}(t)~\cdots~\hat{p}_{M}(t))'$ with
\begin{equation}\label{up01}
\hat{p}_{i}(t)=\begin{cases}
\frac{(t-1)\hat{p}_{i}(t-1)+1}{t},~i=h(t-1), \\
\frac{(t-1)\hat{p}_{i}(t-1)}{t},~~i\neq h(t-1).
\end{cases}
\end{equation}
\end{remark}

Define $\hat{\mu}_w(t)=\mathbf{W}\widehat{\mathbf{p}}_{\min}$. Then, we have
\begin{equation*}\begin{aligned}
{\bf E}[\hat{\mu}_w(t)]=\mu_w,~{\bf E}\big[(\hat{\mu}_w(t)-\mu_w)(\hat{\mu}_w(t)-\mu_w)'\big]=\frac{1}{t}C_w.
\end{aligned}
\end{equation*}
In this case, $\hat{\mu}_w(t)$ is the LMMSUE of $\mu_w$.
Next, based on the LMMSUE, we derive a learning-based $\mathcal {U}_{ad}$ admissible control policy that is optimal for the unknown statistics case.

\begin{theorem}\label{thm02}
For the online quadratic optimization problem with asymmetric information structure, the optimal online control policy in $\mathcal {U}_{ad}$ is designed as
\begin{align}
\hat{u}(t)
=& -\Upsilon_T(t)^{-1} \big( B'P_T(t+1)Ax(t)+B' P_T(t+1)\hat{\mu}_w(t) \nonumber \\
&+B' L_T(t+1)\hat{\mu}_w(t) \big), \label{input02}
\end{align}
where $\widehat{\mathbf{p}}_{\min}(t)$ is the LMMSUE (\ref{est001}) and $\hat{\mu}_w(t)=\mathbf{W}\widehat{\mathbf{p}}_{\min}$.
Moreover, the optimal online index value in (\ref{cost01}) is
\begin{align}
&J_T(\hat{u}(t))\nonumber\\
=&~J_T^*+{\bf Tr} \left({D}_T(0) \mu_w \mu_w'  \right)
+\sum^{T}_{t=1} \frac{1}{t} {\bf Tr}\left( {D}_T(t)C_w \right),\label{cost03}
\end{align}
where $\Upsilon_T(t)$, $P_T(t)$, $L_T(t)$, $H_T(t)$ satisfy the iterative equations (\ref{ric00})-(\ref{ric03}) and
\begin{align}
\mathcal {D}_T(t)=&~(P_T(t+1)+L_T(t+1))'B \Upsilon_T(t)^{-1} B' \nonumber\\
&\times (P_T(t+1)+L_T(t+1))\geq 0. \label{ric04}
\end{align}
\end{theorem}

\textbf{Proof.} See Appendix~\ref{app01b}. \hfill $\Box$

To better understand the performance of the proposed online policy, we need to carry out a detailed regret analysis. For convenience, we state the following hypotheses

{\em
H1) $Q(t)=Q\geq 0$, $R(t)=R>0$ and $P_{T+1}=0$;

H2) $(A,B)$ is stabilizable and $(A,Q^{\frac{1}{2}})$ is observable.}

\begin{lemma}\label{lemma02}
Suppose $P_T(t)$ is the unique positive semi-definite solution to the Riccati equation (\ref{ric01}).
Under hypotheses H1)-H2), $P_T(t)$ is bounded and monotonically nondecreasing as time decreases.
Moreover, when $t\rightarrow -\infty$, $P_T(t)$ converges to the unique solution $\hat{P}>0$ to the following algebraic Riccati equation (ARE)
\begin{equation}\label{are01}
\hat{P}=A'\hat{P}A+Q-A'\hat{P}B(R+B'\hat{P}B)^{-1}B'\hat{P}A.
\end{equation}
\end{lemma}

\textbf{Proof.} See Appendix~\ref{app02}. \hfill $\Box$

\begin{theorem}\label{thm03}
Under hypotheses H1)-H2),  the regret $Reg_T(\hat{u}(t))$ satisfies
\begin{equation}\label{regret01}
Reg_T(\hat{u}(t))\leq O(\ln(T)).
\end{equation}
\end{theorem}

\textbf{Proof.}
By Lemma \ref{lemma02}, $P_T(t)$ is uniformly bounded by $0\leq P_T(t)\leq \hat{P}$, where $\hat{P}$ is the unique positive definite solution satisfying the ARE (\ref{are01}).
By (\ref{ric02b}), $L_T(t)$ is uniquely determined by $P_T(s)$, $s=t+1,\cdots,T$, and thus bounded.
Moreover, by (\ref{ric04}), $\mathcal {D}_T(t)\geq 0$ is determined by $P_T(s)$, $s=t+1,\cdots,T$ and also bounded.
For $\mathcal {D}_T(t)\geq 0$ and $C_w\geq 0$, there exists a constant $\hat{c}>0$ such that
\begin{equation}
{\bf Tr}\left(\mathcal {D}_T(t)C_w \right)\leq \hat{c}.
\end{equation}
By Theorem \ref{thm02}, the regret satisfies
\begin{align}
Reg_T(u(t))
\leq & ~{\bf Tr} \left( {D}_T(0) \mu \mu'  \right) +\sum^{T}_{t=1} \frac{1}{t} \hat{c} \nonumber \\
=& ~{\bf Tr} \left({D}_T(0) \mu \mu'  \right) +\ln(T+r_T)\hat{c}, \label{re01}
\end{align}
where $\lim_{T\rightarrow\infty}r_T=r$ and $r>0$ is the Euler constant.
It follows that $Reg_T(u(t))\leq O(\ln T)$. \hfill $\Box$

\begin{remark}
In our previous work \cite{taniet2017}, a Gittins Index based heuristic policy was developed for a class of pursuit-evasion problems modelled in (\ref{sys00a})-(\ref{cost00}).
The starting point is to minimize an one-step utility function $W(t)$ in each time $t$ as a surrogate cost function, which provides an upper bound of the index function.
Its regret is shown to grow at a linear rate.
Therefore, the proposed learning based control policy in Theorem \ref{thm02} outperforms the Gittins Index based policy in \cite{taniet2017}.
\end{remark}

Suppose $T>0$ is sufficiently large.
Next, we analyse the efficiency of the proposed online control policy $\hat{u}(t)$ compared with the other type of admissible control policies.

{\em Case 1}: Consider the following admissible control policy based on the linear biased estimation defined below
\begin{align}
u_1(t)=& -\Upsilon_T(t)^{-1} ( B'P_T(t+1)Ax(t)+B'(P_T(t+1) \nonumber \\
&+L_T(t+1)) \tilde{\mu}_w(t)),~t=1,2,\cdots,T, \label{case02}
\end{align}
where $\tilde{\mu}_w(t)=\mathbf{W}\tilde{\mathbf{p}}(t)$ and $\tilde{\mathbf{p}}(t)$ is a linear biased estimate satisfying
\be
\tilde{\mathbf{p}}(t)=\sum_{i=0}^{t-1}\tilde{c}_i(t)\xi(i),~t=1,2,\cdots,T.
\ee
In this case, the one-step regret satisfies
\begin{align}
reg_T(t,u_1(t))=&\sum_{i=0}^{t-1} \tilde{c}^2_i(t) {\bf Tr}({D}_T(t)C_w) \nonumber \\ &+(\sum_{i=0}^{t-1} \tilde{c}_i(t)-1)^2 {\bf Tr}\left({D}_T(t)\mu_w\mu_w' \right). \label{estt01}
\end{align}
The minimum regret value of (\ref{estt01}) achieved at
$$\tilde{c}^*_i(t)=\frac{ {\bf Tr} \left({D}_T(t)\mu_w\mu_w' \right)}
{{\bf Tr}\left( t {D}_T(t)\mu_w\mu_w' \right)+{\bf Tr}({D}_T(t)C_w) }.$$
However, the exact values of $\mu_w$ and $C_w$ are unknown to the decision maker, it is infeasible to apply the proposed linear minimum mean square biased estimate (LMMSBE).

{\em Case 2}: Suppose that the decision maker will terminate updating the estimate after a critical time $\bar{t}$, $1\leq \bar{t} <T$.
That is to say, for $0\leq t\leq \bar{t}$, $u_2(t)=\hat{u}(t)$, and for $\bar{t}< t\leq T$,
\begin{align}
{u}_2(t)=&-\Upsilon_T(t)^{-1} ( B'P_T(t+1)Ax(t)+B'P_T(t+1)\hat{\mu}_w(\bar{t}) \nonumber \\
&+B'L_T(t+1)\hat{\mu}_w(\bar{t})). \label{case03}
\end{align}
From the proof of Theorem \ref{thm02}, the regret satisfies
\begin{align*}
Reg_T(u_2(t)) =&~{\bf Tr} \left({D}_T(0) \mu_w \mu_w'  \right)
+\sum^{\bar{t}}_{t=1} \frac{1}{t} {\bf Tr}({D}_T(t)C_w ) \\
&+\sum^{T}_{t=\bar{t}+1} \frac{1}{\bar{t}} {\bf Tr}({D}_T(t)C_w ),
\end{align*}
which implies that $Reg_T(\hat{u}(t))\leq Reg_T(u_2(t)).$
The online control policy $\hat{u}(t)$ in Theorem \ref{thm02} offers a better performance than $u_2(t)$.

{\em Case 3}: Suppose the decision maker does not estimate the probability $\mathbf{p}_w$ and only utilize the state feedback control policy $u(t)=Kx(t)$.
In this case, the optimal feedback control policy is derived as
\begin{equation}\label{case04}
u_3(t)=-\Upsilon_T(t)^{-1} B'P_T(t+1)Ax(t).
\end{equation}
It follows that
\begin{equation*}
{u}_3(t)-u^*(t)=\Upsilon_T(t)^{-1} B' (P_T(t+1)+L_T(t+1))\mu_w ,
\end{equation*}
and
\begin{equation}
reg_T(t,{u}_3(t))
={\bf Tr} \left({D}_T(t) \mu_w \mu_w' \right).
\end{equation}
If $T>0$ is sufficiently large, there exists a critical time $1\leq t_c<T$ such that
\begin{equation}\begin{aligned}
reg_T(t,\hat{u}(t))\leq reg_T(t,{u}_3(t)),~t_c\leq t\leq T.
\end{aligned}\end{equation}
Moreover, under hypotheses H1)-H2), the regret of ${u}_3(t)$ is shown to be linear, which indicates that our policy $\hat{u}(t)$ in Theorem 2 offers a better performance than $u_3(t)$.

\section{Output Feedback Control with Learning}

\subsection{Problem Formulation}

Consider the following discrete time dynamic system
\begin{align}
x(t+1)=&~Ax(t)+Bu(t)+w(t), \label{sys01a}\\
y(t)=&~Cx(t)+v(t),\label{sys01b}
\end{align}
where $y(t) \in\mathbb{R}^n$ is the measurement and $C\in\mathbb{R}^{n\times n}$ is nonsingular with the compatible dimension.
The initial state $x_0\in\mathbb{R}^n$ is a Gaussian random vector with
\begin{align}
\mu_0={\bf E}[w(t)],~C_0= {\bf E} \big[(x_0-\mu_0)(x_0-\mu_0)'\big].
\end{align}
The measurement noise $v(t)$'s are bounded and i.i.d. stochastic process \cite{ghao2001} with
\begin{align}
&\max_{i}\|v(t)\|\leq v_b<\infty,\\
0=&~{\bf E}[v(t)],~Q_v= {\bf E} \big[v(t)v(t)'\big].
\end{align}
We assume that the $w(t)$'s are bounded and form an i.i.d. stochastic process satisfying (\ref{wt01})-(\ref{wt02}).
The random variables $x_0$, $w(t)$, $v(t)$ are assumed to be mutually independent.
Moreover, we emphasize that the probability $\mathbf{p}_w$ is a priori unknown to the decision maker.
The objective is to minimize the index function (\ref{cost01}) with asymmetric information structure.

Generally speaking, to solve the quadratic optimization problem (\ref{cost01}) subject to (\ref{sys01a})-(\ref{sys01b}), one could apply the well known Kalman filter to estimate the value of the state $x(t)$ and based on that design the optimal offline control policy to minimize the index function.
To be specific, denote $Y(t)$ to be the observation set $\{y(0),y(0),\cdots,y(t)\}$.
Define $\hat{x}_{t|t-1}={\bf E}[x(t)|Y(t-1)]$, $\hat{x}_{t|t}={\bf E}[x(t)|Y(t)]$ and
\begin{align}
\Lambda_{t|t-1}=&~{\bf E}[(x(t)-\hat{x}_{t|t-1})(x(t)-\hat{x}_{t|t-1})'|Y(t-1)], \\
\Lambda_{t|t}=&~{\bf E}[(x(t)-\hat{x}_{t|t})(x(t)-\hat{x}_{t|t})'|Y(t)],
\end{align}
where ${\bf E} [x|Y]$ defines the conditional expectation of the random variable $x$ w.r.t. $Y$.
Applying the standard Kalman filtering \cite{kalman} yields that
\begin{align}
\hat{x}_{t|t}=&~\hat{x}_{t|t-1}+\Lambda_{t|t-1}C'(C\Lambda_{t|t-1}C'+Q_v)^{-1} \nonumber \\
&\times (y(t)-C\hat{x}_{t|t-1}), \\
\hat{x}_{t+1|t}=&~ A\hat{x}_{t|t}+Bu(t)+\mu_w,  \label{est01}
\end{align}
where
\begin{align}
\Lambda_{t|t}=& ~\Lambda_{t|t-1}-\Lambda_{t|t-1}C'(C\Lambda_{t|t-1}C'+Q_v)^{-1} C\Lambda_{t|t-1},\nonumber \\
\Lambda_{t+1|t}=&~ A\Lambda_{t|t}A'+C_w.
\end{align}
The initial conditions are
\be\hat{x}_{0|-1}={\bf E}[x_0]=\mu_0,~\Lambda_{0|-1}={\bf E}[(x_0-\mu_0)(x_0-\mu_0)'].\ee
By utilizing the separation principle, the optimal offline control policy is derived as
\begin{align}
u^*(t)=&-\Upsilon_T(t)^{-1} \big( B'P_T(t+1)A \hat{x}_{t|t}+B'P_T(t+1)\mu_w \nonumber \\
&+B'L_T(t+1)\mu_w \big).
\end{align}

In the current model, since the exact values of $\mu_w$ and $C_w$ are unknown, the classic Kalman filter and the optimal offline control strategy cannot be applied for the asymmetric information case.
Instead, we introduce an one-step state estimation based on the observation $y(t)$ at each time $t=1,2,\cdots,T$.
The original problem is reduced to a quadratic optimization problem with a non-white system noise \cite{color1995}.
This modified optimization problem is challenging.
In this study, we derive a suboptimal offline control policy conditioned on the assumption that the one-step state estimation is applied and the probability statistics of the system are known.
Based on the LMMSUE, we propose a learning based online control policy.
The quasi-regret between the online known statistics cost and the heuristic offline unknown statistics suboptimal cost is shown to be bounded by $O(\ln T)$.

\subsection{Learning Based Control Policy and Regret Analysis}
With the output dynamic equation (\ref{sys01b}), we introduce a simple one-step state estimate
\begin{equation}
\hat{x}(t)={\bf E}[x(t)|y(t)]=C^{-1}y(t),~t=0,1,\cdots,T,
\end{equation}
which implies that
\begin{align}
\hat{x}(t+1)=&~A\hat{x}(t)+Bu(t)+s(t),\label{sys02a} \\
s(t)=&~w(t)-AC^{-1}v(t)+C^{-1}v(t+1). \label{sys02b}
\end{align}
In this case, $s(t)$ is a colored noise with $\mu_s={\bf E}[s(t)]=\mu_w$ and
\begin{align}
Q_s={\bf E}[s(t)s'(t)]=Q_w+C^{-1}Q_vC^{-1}+AC^{-1}Q_vC^{-1}A'.
\nonumber
\end{align}
Moreover, the error covariance is
\begin{equation*}
\Lambda(t)={\bf E} \left[(x(t)-\hat{x}(t))(x(t)-\hat{x}(t))'\right]=C^{-1}Q_v C^{-1}.
\end{equation*}
The index function (\ref{cost01}) can be rewritten as
\begin{align}
J_T(u(t))= &\sum_{t=0}^{T} {\bf E} \big[ \hat{x}'(t)Q(t)\hat{x}(t)+u'(t)R(t)u(t) \big] \nonumber \\
&+{\bf E} \big[ \hat{x}'(T+1)P_{T+1}\hat{x}(T+1)\big]
-D_T, \label{cost01a}
\end{align}
where
\begin{align}
D_T=&~\sum_{t=0}^{T}{\bf Tr}(Q(t)\bar{Q}_v)+{\bf Tr}(P_{T+1}\bar{Q}_v), \\
\bar{Q}_v=&~C^{-1}Q_vC^{-1}.
\end{align}
The original quadratic optimization problem (\ref{cost01}) is reduced to minimizing (\ref{cost01a}) w.r.t. (\ref{sys02a})-(\ref{sys02b}).
By utilizing the one-step state estimate, we derive a heuristic suboptimal offline result.

\begin{theorem}
\label{thm04}
Suppose the probability $\mathbf{p}_w$ is known a priori.
A suboptimal offline control policy of the quadratic optimization problem (\ref{cost01}) is given by
\begin{align}
u_a(t)=&-\Upsilon_T(t)^{-1} \big( B'P_T(t+1)A \hat{x}(t) +B'P_T(t+1)\mu_w \nonumber \\
&+B'L_T(t+1)\mu_w \big), \label{in02}
\end{align}
while the index value of (\ref{cost01}) is
\begin{equation}
J_T(u_a(t))= \hat{x}'(0)P_T(t)\hat{x}(0)+2\hat{x}'(0)L_T(0)\mu_w +H_T,
\end{equation}
where $\Upsilon_T(t)$, $P_T(t)$, $L_T(t)$ satisfy (\ref{ric00})-(\ref{ric02}) and
\begin{align}
H_T=&\sum_{t=0}^{T}\Big\{ -\mu_w'(P_T(t+1)+L_T(t+1))'B \Upsilon_T(t)^{-1} B' \nonumber\\
&\times (P_T(t+1)+L_T(t+1))\mu_w+2\mu_w' L_T(t+1)\mu_w \nonumber \\
&+{\bf Tr} \left( A'P_{T}(t+1)B \Upsilon_T(t)^{-1} B'P_{T}(t+1)A \bar{Q}_v \right) \nonumber \\
&+{\bf Tr} \left(P_T(t+1)Q_w\right)\Big\}
-{\bf Tr}(P_T(0)\bar{Q}_v).\label{ric05}
\end{align}
\end{theorem}

\textbf{Proof.} See Appendix~\ref{app03}. \hfill $\Box$

Next, we study the LMMSUE of $\mathbf{p}_w$.
Due to the presence of the measurement noise $v(t)$, at each time $t=1,2,\cdots, T$, it is difficult to reach the exact value of $w(t-1)$.
To guarantee the exact observation of $w(t-1)$, we state the following hypothesis

{\em H3) For each $i,j=1,2,\cdots,M$ and $i\neq j$,
\begin{equation}
\|w_i-w_j\|> \frac{2(1+\|A\|)v_b}{\|C\|}.
\end{equation}
}

At each time $t=1,2,\cdots,T$, define $\hat{w}(t-1)=\hat{x}(t)-A\hat{x}(t-1)-Bu(t-1)$.
It follows that
\begin{align}
\|\hat{w}(t-1)-w(t-1) \|=&~ \| C^{-1}v(t)-AC^{-1}v(t-1) \| \nonumber\\
 \leq &~  \frac{(1+\|A\|)v_b}{\|C\|}.
\end{align}
Suppose that $w(t-1)=w_{h(t-1)}$. For each $i\neq h(t-1)$, we obtain
\begin{align*}
\frac{2(1+\|A\|)v_b}{\|C\|}<& \|{w}_i-w_h \| \\
\leq & \|{w}_i-\hat{w}(t-1) \|+\|w_h-\hat{w}(t-1) \| \\
\leq & \|{w}_i-\hat{w}(t-1) \| +\frac{(1+\|A\|)v_b}{\|C\|},
\end{align*}
which implies that
\begin{equation}
\|\hat{w}(t-1)-w_i \| >  \frac{(1+\|A\|)v_b}{\|C\|}\geq  \| \hat{w}(t-1) -w_h \|.
\end{equation}
Therefore, at each time $t=0,1,\cdots,T-1$, we have $w(t)=w_{h(t)}$, where $h(t)$ is determined by
\begin{equation}
h(t)={\arg\min}_{i=1,\cdots,M}\|\hat{w}(t)-w_i \|.
\end{equation}
By Remark \ref{remark02}, we have $\xi(t)=\xi_{h(t)}$ and update $\widehat{\mathbf{p}}_{\min}(t+1)=[\hat{p}_{1}(t+1)~\hat{p}_{2}(t+1)~\cdots~\hat{p}_{M}(t+1)]'$ with
\begin{equation}\label{up01}
\hat{p}_{i}(t+1)=\begin{cases}
\frac{t\hat{p}_{i}(t)+1}{t+1},~i=h(t), \\
\frac{t\hat{p}_{i}(t)}{t+1},~~i\neq h(t).
\end{cases}
\end{equation}

Based on the LMMSUE, we are in a position to present a learning based online control policy as follows.
\begin{theorem}\label{thm05}
Suppose the probability $\mathbf{p}_w$ is unknown.
Under hypotheses H3), an admissible online control policy is derived as
\begin{align}
\hat{u}_a(t)
=& -\Upsilon_T(t)^{-1} \big( B'P_T(t+1)A\hat{x}(t)+B'P_T(t+1)\hat{\mu}_w(t) \nonumber \\
&+B'L_T(t+1) \hat{\mu}_w(t) \big), \label{in01}
\end{align}
while the index value in (\ref{cost01}) is
\begin{align}
&J_T(\hat{u}_a(t))
=\hat{x}'(0)P_T(t)\hat{x}(0)+2\hat{x}'(0)L_T(0)\mu_w+H_T \nonumber\\
&~~~~~~~~~~+{\bf Tr} \left({D}_T(0) \mu_w \mu_w'  \right)
+\sum^{T}_{t=1} \frac{1}{t} {\bf Tr}\left( {D}_T(t)C_w \right), \label{co01}
\end{align}
where $\Upsilon_T(t)$, $P_T(t)$, $L_T(t)$ satisfy (\ref{ric00})-(\ref{ric01}), $\mathcal {D}_T(t)$ satisfies (\ref{ric04}) and $H_T$ satisfies (\ref{ric05}).
Moreover, under hypotheses H1)-H2), the quasi-regret between the online cost $J_T(\hat{u}_a(t))$ and the offline cost $J_T(u_a(t))$ satisfies
\begin{align}
\bar{R}eg_T(\hat{u}_a(t))\leq O(\ln T).
\end{align}
\end{theorem}

\textbf{Proof.} See Appendix~\ref{app04}. \hfill $\Box$

\begin{remark}
In this paper, for the sake of simplicity, we investigate the single-armed optimization problem with asymmetric information structure \cite{taniet2017}.
From a general application perspective, it is of interest to consider the following model
\begin{align}
x(t+1)=&~Ax(t)+u(i(t),t)+w_{i(t)}(t), \label{sys10a}\\
y(t)=&~Cx(t)+v(t),\label{sys10b}
\end{align}
where $i(t)\in\mathbb{K} \triangleq\{1,2,\cdots,K\}$ is the control mode.
Moreover, the system matrix $A$ and the index matrices $Q\geq 0$, $R>0$ are assumed to be unknown.
The first predator-prey example is a special case of (\ref{sys10a})-(\ref{sys10b}).
How to best utilize the observed trajectories to estimate the unknown information and based on that propose a learning based control policy is a challenging future work direction.
\end{remark}

\section{Illustrative Examples}

In this section, we present two numerical examples to illustrate the effectiveness of our theoretical results.

\begin{example}
Consider the predator-prey model in (\ref{sys00a})-(\ref{cost00}).
For convenience, we simply set $\beta=1$ and $\mathbb{R}^n=\mathbb{R}^2$.
Assume the initial positions are $x_p=[1~0]'$ and $x_e=[0~0]'$.
The prey has the following four evading policies $$v_1=[1~0]',~v_2=[-1~0]',~v_3=[0~1]',~v_4=[0~-1]'.$$
with the evading probability distribution
$$\mathbf{p}_v=[0.2~0.1~0.6~0.1]'.$$
In this case, we obtain
\[ \begin{matrix}
\mu_v=\begin{bmatrix}
 0.1 \\ 0.5
\end{bmatrix},
~Q_v=\begin{bmatrix}
 0.3 & 0 \\ 0& 0.7
\end{bmatrix}.
\end{matrix}\]
Note that the evading probability distribution $\mathbf{p}_v$ is unknown to the predator.

If we set $x(t)=x_p(t)-x_e(t)$ and $T=200$, the first predator-prey problem in (\ref{sys00a})-(\ref{cost00}) can be reformulated as the state feedback case (\ref{cost01}) with
\[ \begin{matrix} A=B=\begin{bmatrix}
 1 & 0 \\ 0& 1
\end{bmatrix},
~Q=R=\begin{bmatrix}
 1 & 0 \\ 0& 1
\end{bmatrix},~
x(0)=\begin{bmatrix}
 1 \\ 0
\end{bmatrix}.
\end{matrix}\]
It follows that  $(A,B)$ is stabilizable and $(A,Q^{\frac{1}{2}})$ is observable.
Suppose the evading probability distribution $\mathbf{p}_v$ is known to the predator. By Theorem \ref{thm01}, we obtain the optimal offline control policy $u^*(t)$ in (\ref{input}) which minimizes the index function (\ref{cost01}) with $J_T^*=292.1660$.

By utilizing the proposed admissible control policy $\hat{u}(t)$ in (\ref{input02}) with the LMMSUE $\widehat{\mathbf{p}}_{\min}(t)$, we obtain the index cost with $J_T(\hat{u}(t))=304.2107$.
Thus, the regret is $$Reg_T(\hat{u}(t))=J_T(\hat{u}(t))-J_T^*=12.0446.$$
We propose the trajectories of the one-step regret $reg_T(\hat{u(t)})$ as shown in Fig. 2,
where $reg_T(\hat{u(T)})=0$ due to the terminal conditions $P_{T+1}=L_{T+1}=0$.
\begin{figure}[thpb]
      \centering
      \includegraphics[scale=0.5]{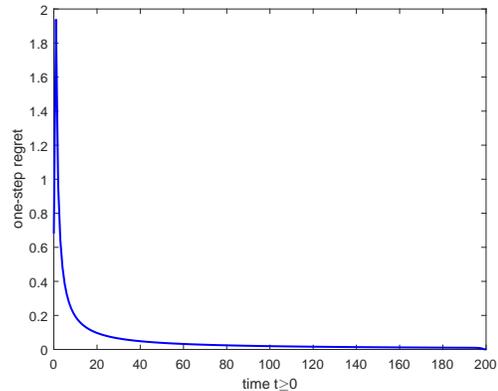}
      \caption{Trajectories of one-step regret}
\end{figure}

Moreover, define the regret percentage to be
\begin{align}
c_T(\hat{u}(t))=\frac{Reg_T(\hat{u}(t))}{T}\times 100\% .
\end{align}
For different terminal time $T>0$, the optimal offline index value $J_T^*$, the optimal online index value $J_T(\hat{u}(t))$, the regret $Reg_T(\hat{u}(t))$ and the percentage $c_T(\hat{u}(t))$ can be summarized in Table I.
It can be concluded that the regret of the proposed online control policy grows at a sub-linear rate.
\begin{table}[thpb]
\centering
\caption{Regret Analysis in Example 1}
  \begin{tabular}{|c|c|c|c|c|}
  \hline
    & ${J}^*_T$ & ${J}_T(\hat{u})$ & $Reg_T(\hat{u})$ & $c_T(\hat{u})$\\
    \hline
  \hspace{-1mm}$T=20$ & \hspace{-1mm}$29.8439$ & \hspace{-1mm}$37.2439$ & \hspace{-1mm}$7.4000$ & \hspace{-1mm}$37.0000\%$ \\
  \hline
  \hspace{-1mm}$T=50$ & \hspace{-1mm}$73.5643$ & \hspace{-1mm}$82.8646$ & \hspace{-1mm}$9.3004$ & \hspace{-1mm}$18.6008\%$ \\
  \hline
  \hspace{-1mm}$T=100$ & \hspace{-1mm}$146.4315$ & \hspace{-1mm}$157.1141$ & \hspace{-1mm}$10.6825$ & \hspace{-1mm}$10.6825\%$  \\
  \hline
  \hspace{-1mm}$T=200$ & \hspace{-1mm}$292.1660$ & \hspace{-1mm}$304.2107$ & \hspace{-1mm}$12.0446$ & \hspace{-1mm}$6.0223\%$  \\
  \hline
  \hspace{-1mm}$T=500$ & \hspace{-1mm}$729.3696$ & \hspace{-1mm}$743.2008$ & \hspace{-1mm}$13.8312$ & \hspace{-1mm}$2.7662\%$  \\
  \hline
  $T=1000$ & \hspace{-1mm}$1458.0422$ & \hspace{-1mm}$1473.2200$ & \hspace{-1mm}$15.1778$ & \hspace{-1mm}$1.5178\%$  \\
\hline
  \hspace{-1mm} $T=2000$ & \hspace{-1mm}$2915.3873$ & \hspace{-1mm}$2931.9100$ & \hspace{-1mm}$16.5227$ & \hspace{-1mm} $0.8261\%$  \\
\hline
\end{tabular}
\end{table}
\end{example}

\begin{example}
Consider the modified product pricing model in (\ref{dy01})-(\ref{dy02}).
We assume that $\alpha=\frac{1}{4}$, $b=2$ and $e(t)$'s are bounded and i.i.d. stochastic process with
\begin{align*}
&e_1=0,~e_2=0.1,~e_3=-0.1,~e_4=0.2,~e_5=-0.2,\\
&~~~~e_6=0.3,~e_7=-0.3,~e_8=0.4,~e_9=-0.4.
\end{align*}
The probability distribution $\mathbf{p}_e$ is assume to be
$$\mathbf{p}_e=[0.25~0.15~0.15~0.1~0.1~0.075~0.075~0.05~0.05]'.$$
Moreover, we assume that $c_1=c_2=c_3=1$.

If we set $X(t)=[y(t)~z(t)]'$, $U(t)=[v(t)~u(t)]'$ and $W(t)=[w(t)~0]'$, the second product pricing problem can be reformulated as the state feedback case (\ref{cost01}) with
\[ \begin{matrix} A=\begin{bmatrix}
 0 & 0 \\ 0& 1
\end{bmatrix},
~B=\begin{bmatrix}
 1 & 0 \\ 0& 1
\end{bmatrix},
\end{matrix}\]
\[ \begin{matrix}
Q=\begin{bmatrix}
 \frac{1}{16} & -\frac{1}{4} \\ -\frac{1}{4} & 1
\end{bmatrix}\geq 0,
~R=\begin{bmatrix}
 1 & 0 \\ 0 & 1
\end{bmatrix}>0.
\end{matrix}\]
It follows that $(A,B)$ is stabilizable and $(A,Q^{\frac{1}{2}})$ is observable.
Moreover, the mean and variance of $W(t)$ is
\[ \begin{matrix}
\mu_W=\begin{bmatrix}
 3.6 \\ 0
\end{bmatrix},
~Q_W=\begin{bmatrix}
 13.6080 & 0 \\ 0& 0
\end{bmatrix}.
\end{matrix}\]

For different terminal time $T>0$, it follows from Theorem \ref{thm01}-\ref{thm03} that the optimal index value $\bar{J}_T^*$, the index value $\bar{J}_T(\hat{u}(t))$, the regret $Reg_T(\hat{u}(t))$ and the percentage $c_T(\hat{u}(t))$ can be summarized in Table II.

\begin{table}[thpb]
\centering
\caption{Regret Analysis in Example II}
  \begin{tabular}{|c|c|c|c|c|}
  \hline
 & $\bar{J}^*_T$ & $\bar{J}_T(\hat{u})$ & $Reg_T(\hat{u})$ &  $c_T(\hat{u})$\\
    \hline
  $T=20$ & $1.3786$ & $2.3309$ & $0.9524$ & $4.7620\%$ \\
  \hline
  $T=50$ & $2.5936$ & $3.5845$ & $0.9909$ & $1.9818\%$ \\
  \hline
  $T=100$ & $4.6186$ & $5.6380$ & $1.0194$ & $1.0194\%$ \\
  \hline
  $T=200$ & $8.6686$ & $9.7163$ & $1.0477$ & $0.5239\%$\\
  \hline
  $T=500$ & $20.8186$ & $21.9036$ & $1.0850$ & $0.2170\%$ \\
  \hline
  $T=1000$ & $41.0686$ & $42.1817$ & $1.1131$ & $0.1113\%$ \\
\hline
  $T=2000$ & $81.5686$ & $82.7098$ & $1.1412$ & $0.0571\%$\\
\hline
\end{tabular}
\end{table}
\end{example}

\section{Conclusions}

In this paper, we focus on an online quadratic optimization problem with asymmetric information structure.
We assume that a single predator with rich information input is pitted against a single prey with limited input information.
Motivated by the OCO methodology, we develop a robust approach that enable the predictor-agent learn the probability statistics of the system with the LMMSUE.
Based on that we propose a learning based optimal online control policy.
Its regret grows at a sub-linear rate, and is shown to be bounded by $O(\ln T)$,
which implies the online performance can converge asymptotically to that of the offline optimal performance.

As future work, there are two promising research directions.
The first research direction is to figure out more optimal online control strategies and analysis framework for existing online quadratic optimization problems.
The other direction is to extend the two-player models to more complicated models such as multi-agent systems.
With unknown statistics of multiplicative noise or network topology, it is infeasible to utilize the classic distributed control strategies.
The online optimization approach can offer a promising but challenging new direction.

\appendices

\section{Proof of Theorem \ref{thm01}}\label{app01}
\textbf{Proof.}
The proof is based on the dynamic programming approach.
For each time $t=0,1,\cdots,T$, define the following cost-to-go function
\begin{equation}\label{costgo01}
\begin{aligned}
\mathcal{G}(t)=& \min_{u(t)}  G(t),
\end{aligned}\end{equation}
where
\begin{equation}\label{costgo01a}
G(t)= {\bf E} \big[ x'(t)Q(t)x(t)+u'(t)R(t)u(t)+\mathcal{G}(t+1)\big],
\end{equation}
and
\begin{equation}
\mathcal{G}(T+1)={\bf E}\left[ x'(T+1)P_{T+1}x(T+1) \right].
\end{equation}
Next, we show that
\begin{align}
 \mathcal{G} (t)=&~{\bf E}\left[ x'(t)P_T(t)x(t)+2x'(t)L_T(t)\mu_w \right] \nonumber \\
&+\sum_{j=t}^TH_T(j),
\end{align}
where $P_T(t)$, $L_T(t)$ and $H_T(t)$ satisfy (\ref{ric01})-(\ref{ric03}).

For $t=T$, it follows from (\ref{costgo01a}) that
\begin{equation*}
\begin{aligned}
&G(T) \\
=&~{\bf E} \Big[ x'(T)Q(T)x(T)\hspace{-0.5mm}+\hspace{-0.5mm}u'(T)R(T)u(T)
+(Ax(T)\hspace{-0.5mm}+\hspace{-0.5mm} Bu(T) \\
&+w(T))'P_{T+1} (Ax(T)+Bu(T)+w(T))\Big] \\
=&~{\bf E} \Big[ (u(T)\hspace{-0.6mm}+\hspace{-0.6mm}\Upsilon_T(T)^{-1}
\hspace{-0.35mm}(B'P_{T+1}Ax(T)\hspace{-0.6mm}+\hspace{-0.6mm}
B'P_{T+1}\mu_w))'\Upsilon_T(T) \\
&\times(u(T)+\Upsilon_T(T)^{-1}(B'P_{T+1}Ax(T)+B'P_{T+1}\mu_w)) \\
&-(B'P_{T+1}Ax(T)+B'P_{T+1}\mu_w)'\Upsilon_T(T)^{-1}\\
&\times(B'P_{T+1}Ax(T)+B'P_{T+1}\mu_w)  \\
&+x'(T)(Q(T)+A'P_{T+1}A)x(T) + 2x'(T)A'P_{T+1}\mu_w \Big] \\
&+{\bf Tr}(P_{T+1}Q_w),
\end{aligned}\end{equation*}
where $\Upsilon_T(T)=R(T)+B'P_{T+1}B$. At time $T$, the optimal control policy $u^*(T)$ is
\begin{equation}\label{con11}
u^*(T)=-\Upsilon_T(T)^{-1}(B'P_{T+1}Ax(T)+B'P_{T+1}\mu_w),
\end{equation}
while the cost-to-go function $\mathcal{G}(T)$ satisfies

\begin{equation*}\label{eq11a}
\begin{aligned}
& \mathcal{G}(T)=\min_{u(T)} G(T)\\
=& ~{\bf E} \Big[x'(T)(Q(T)+A'P_{T+1}A- A'P_{T+1}B \Upsilon_T(T)^{-1}B'  \\
& \times P_{T+1}A)x(T)+ 2x'(T)(A'-A'P_{T+1}B \Upsilon_T(T)^{-1} B') \\
&\times P_{T+1}\mu_w \Big]
-\mu_w' P_{T+1}B\Upsilon_T(T)^{-1}B'P_{T+1}\mu_w \\
&+{\bf Tr}(P_{T+1}Q_w) \\
=&~{\bf E}\big[ x'(T)P_T(T)x(T)+2x'(T)L_T(T)\mu_w \big]+H_T(T),
\end{aligned}
\end{equation*}
where the parameters $P_T(T)$, $L_T(T)$ and $H_T(T)$ satisfy (\ref{ric01})-(\ref{ric03}) with $t=T$.

For each $t=0,1,\cdots,T-1$, suppose
\begin{align*}
 \mathcal{G} (t+1)  =&~{\bf E}[ x'(t+1)P_T(t+1)x(t+1) \nonumber \\
&+2x'(t+1)L_T(t+1)\mu_w ] +\sum_{j=t+1}^TH_T(j).
\end{align*}
It follows that
\begin{align*}
&G(t)\\
=& ~{\bf E}\Big[x'(t)Q(t)x(t)+u'(t)R(t)u(t) +x'(t+1) P_T(t+1)\\
&\times x(t+1)+2x'(t+1)L_T(t+1)\mu_w \Big]+\sum_{j=t+1}^TH_T(j)\\
=&~{\bf E} \Big[ (u(t)+\Upsilon_T(t)^{-1}(B'P_T(t+1)Ax(T)+B'(P_T(t+1) \\
&+L_T(t+1))\mu_w))'\Upsilon_T(t)(u(t)+\Upsilon_T(t)^{-1} \\
&\times(B'P_T(t+1)Ax(t)
\hspace{-0.5mm}+\hspace{-0.5mm}B'(P_T(t+1)
\hspace{-0.5mm}+\hspace{-0.5mm}L_T(t+1))\mu_w)) \\
&-(B'P_T(t+1)Ax(t)+B'(P_T(t+1)\hspace{-0.5mm}+\hspace{-0.5mm}L_T(t+1))\mu_w)' \\
&\times \Upsilon_T(t)^{-1}(B'P_T(t+1)Ax(t)+B'(P_T(t+1)\\
&+L_T(t+1))\mu_w) +x'(t)(Q(t)+A'P_T(t+1)A)x(t) \\
&+ 2 x'(t)A'P_T(t+1)\mu _w +2x'(t)A'L_T(t+1)\mu_w \Big]\\
&+2\mu_w'L_T(t+1)\mu_w
+{\bf Tr}(P_T(t+1)Q_w)+\sum_{j=t+1}^TH_T(j).
\end{align*}
Then, the optimal control policy is $u^*(t)$ in (\ref{input}),
which implies that
\begin{equation*}
\begin{aligned}
&\mathcal {G}(t)=\min_{u(t)}G(t) \\
=&~ {\bf E} \Big[x'(t)(Q(t)+A'P_T(t+1)A
\hspace{-0.5mm}-\hspace{-0.5mm} A'P_T(t+1)B \Upsilon_T(t)^{-1} \\
&\times B'P_T(t+1)A)x(t) + 2x'(t)(A'-A'P_T(t+1)B \\
&\times\Upsilon_T(t)^{-1} B')(P_T(t+1)+L_T(t+1))\mu_w \Big] \\
&-\mu_w'(P_T(t+1)
\hspace{-0.5mm}+\hspace{-0.5mm}L_T(t+1))'B\Upsilon_T(t)^{-1}B'(P_T(t+1)\\
&+L_T(t+1))\mu_w +2\mu'L_T(t+1)\mu_w \\
&+{\bf Tr}(P_T(t+1)Q_w)
+\sum_{j=t+1}^TH_T(j)  \\
=&~ {\bf E}\big[ x'(t)P_T(t)x(t)+2x'(t)L_T(t)\mu_w \big]+\sum_{j=t}^TH_T(j).
\end{aligned}
\end{equation*}
Utilize the dynamic programming with the cost-to-go function (\ref{costgo01}) yields the optimal index value satisfies (\ref{cost02}), which completes this proof.
\hfill $\Box$

\section{Proof of Proposition \ref{prop01}}\label{app01a}

\textbf{Proof.}
For each time $t=0,1,\cdots, T$, define the following Lyapunov function
\begin{equation*}\label{vt01}
V(t)={\bf E} \big[ x'(t)P_T(t)x(t)+2x'(t)L_T(t)\mu_w \big]
+\sum_{i=t}^TH_T(i).
\end{equation*}
It follows that
\begin{equation*}
\label{vt02}\begin{aligned}
&V(t)-V(t+1) \\
=&~ {\bf E}\big[ x'(t)P_T(t)x(t)+2x'(t)L_T(t)\mu_w \big]+\sum_{i=t}^T H_T(i) \\
&-{\bf E}\big[ (Ax(t) \hspace{-0.7mm}+\hspace{-0.7mm}Bu(t) \hspace{-0.7mm}+\hspace{-0.7mm}w(t))'
P_T(t\hspace{-0.7mm}+\hspace{-0.7mm}1)
(Ax(t)\hspace{-0.7mm}+\hspace{-0.7mm}Bu(t)\hspace{-0.7mm}+\hspace{-0.7mm}w(t))\\
&+2(Ax(t)+Bu(t)+w(t))'L_T(t+1)\mu_w \big]
\hspace{-0.6mm}-\hspace{-0.9mm}\sum_{i=t+1}^TH_T(i) \\
=& ~ {\bf E} \big[ x'(t)Q(t)x(t)+u'(t)R(t)u(t) \big] \\
&-{\bf E}\big[(u(t)-u^*(t))'\Upsilon_T(t)(u(t)-u^*(t))\big],
\end{aligned}
\end{equation*}
which implies that
\begin{equation*}\label{vt03}\begin{aligned}
&V(0)-V(T+1)=\sum^T_{t=0}V(t)-V(t+1) \\
=& ~x_0'P_T(0)x_0+2x_0'L_T(0)\mu+\sum_{t=0}^T H_T(t) \\
& -{\bf E}[ x'(T+1)P_T(T+1)x(T+1)]\\
=& \sum^T_{t=0} {\bf E}[ x'(t)Q(t)x(t)+u'(t)R(t)u(t)]\\
&-\sum^T_{t=0} {\bf E}[(u(t)-u^*(t))'\Upsilon_T(t) (u(t)-u^*(t))].
\end{aligned}\end{equation*}
In this case, we obtain
\begin{equation*}\label{eq011}\begin{aligned}
J_T(u(t))
=&~ x_0'P_T(0)x_0+2x_0'L_T(0)\mu+\sum_{t=0}^T H_T(t) \\
&+\sum_{t=0}^T {\bf E}\big[ (u(t)-u^*(t))'\Upsilon_T(t) (u(t)-u^*(t)) \big],
\end{aligned}
\end{equation*}
which completes this proof.
\hfill $\Box$

\section{Proof of Theorem \ref{thm02}}\label{app01b}
\textbf{Proof.}
By Proposition \ref{prop01}, we first show that the regret for the developed control policy $\hat{u}(t)$ satisfies
\begin{align}
Reg_T(\hat{u}(t))
=&~{\bf Tr} \left( {D}_T(0) \mu_w \mu_w'  \right)
+\sum^{T}_{t=1} \frac{1}{t} {\bf Tr}\left( {D}_T(t)C_w \right). \label{error01a}
\end{align}
For $t_0=0$, the decision maker has no observation.
With the initial estimate $\widehat{\mathbf{p}}_{\min}(0)=[0~0~\cdots~0]'$, the control policy is designed to be
\begin{equation}\label{input01}
\hat{u}(0)=-\Upsilon_T(0)^{-1} B'P_T(1)Ax(0),
\end{equation}
which is the feedback of the initial state $x(0)=x_0$.
In this case, we have
\begin{equation}
\hat{u}(0)-u^*(0)=\Upsilon_T(0)^{-1} B' (P_T(1)+L_T(1))\mu_w ,
\end{equation}
which implies that
\begin{align}
reg_T(0,\hat{u}(0))
=& ~{\bf E} \left[ (\hat{u}(0)-u^*(0))'\Upsilon_T(0) (\hat{u}(0)-u^*(0)) \right] \nonumber \\
=& ~{\bf Tr} \left( {D}_T(0) \mu_w \mu_w'  \right) ,
\end{align}
where $\mathcal {D}_T(t)\geq 0$ is given in (\ref{ric04}).
For each time $t=1,2,\cdots,T$, the decision maker observes the exact value of $\xi(i)$, $i=0,1,\cdots, t-1$.
In this case,
we obtain
\begin{equation*}
\begin{aligned}
&\hat{u}(t)-u^*(t) \\
=& -\Upsilon_T(t)^{-1}B'(P_T(t+1)+L_T(t+1)) (\hat{\mu}_w(t)-\mu_w),
\end{aligned}\end{equation*}
which implies that
\begin{align}
reg_T(t,\hat{u}(t))=&~{\bf E}\big[ (\hat{u}(t)-u^*(t))'\Upsilon_T(t)(\hat{u}(t)-u^*(t)) \big] \nonumber \\
=&~ {\bf E}\big[(\hat{\mu}_w(t)-\mu_w)'\mathcal {D}_T(t) (\hat{\mu}_w(t)-\mu_w)\big] \nonumber \\
=& ~\frac{1}{t} {\bf Tr}\left( {D}_T(t)C_w \right).
\end{align}
With the updated estimate $\widehat{\mathbf{p}}_{\min}(t)$ and the control policy $\hat{u}(t)$, the regret satisfies (\ref{error01a}).
It follows from Proposition \ref{prop01} that the index value in (\ref{cost01}) is
\begin{align*}
J_T(\hat{u}(t))
=&~x_0'P_T(0)x_0 + 2 x_0'L_T(0) \mu_w+\sum_{t=0}^T H_T(t)\\
&+{\bf Tr} \left({D}_T(0) \mu_w \mu_w'  \right)
+\sum^{T}_{t=1} \frac{1}{t}{\bf Tr}\left(  {D}_T(t)C_w \right).
\end{align*}

For each online control policy $u_1(t)\in \mathcal {U}_{ad}$ satisfying
\begin{align}
u_1(t)=& -\Upsilon_T(t)^{-1} ( B'P_T(t+1)Ax(t)+B'(P_T(t+1) \nonumber \\
&+L_T(t+1)) \check{\mu}_w(t)),~t=1,2,\cdots,T,\label{case01}
\end{align}
where $\check{\mu}_w(t)=\mathbf{W}\check{\mathbf{p}}(t)$ and $\check{\mathbf{p}}(t)$ is a linear unbiased estimate satisfying
$$\check{\mathbf{p}}(t)=\sum_{i=0}^{t-1}\check{c}_i(t)\xi(i),~t=1,2,\cdots,T,$$
with
$\check{c}_i(t)\neq\frac{1}{t}$, $i=0,1,\cdots,t-1$, and $\sum_{i=0}^{t-1}\check{c}_i(t)=1$.
In this case, the regret of $u_1(t)$ satisfies
\begin{equation*}\begin{aligned}
&Reg_T({u}_1(t))\\
=&\sum_{t=0}^T {\bf E}\big[ ({u}_1(t)-u^*(t))'\Upsilon_T(t)({u}_1(t)-u^*(t)) \big] \\
=& ~{\bf Tr} \left({D}_T(0) \mu_w \mu_w'  \right)
\hspace{-0.5mm}+\hspace{-0.65mm}\sum_{t=1}^T \sum_{i=0}^{t-1}\check{c}^2_i(t) {\bf Tr}\left({D}_T(t)C_w \right).
\end{aligned}\end{equation*}
Define
\begin{equation}\label{fc}
f(\check{c})=\sum_{i=0}^{t-1}\check{c}_i^2(t).
\end{equation}
Applying $\check{c}_0(t)=1-\sum_{i=1}^{t-1}\check{c}_i(t)$ to (\ref{fc}), we obtain
\begin{equation}\label{fc01}
f(\check{c})=\sum_{i=1}^{t-1}\check{c}^2_i(t)+(1-\sum_{i=1}^{t-1}\check{c}_i(t))^2.
\end{equation}
For each $j=1,2,\cdots,t-1$, we have
\begin{equation}\label{fc02}
f_{\check{c}_j(t)}(\check{c})=2(\check{c}_j(t)+\sum_{i=1}^{t-1}\check{c}_i(t)-1).
\end{equation}
Suppose $f_{\check{c}_j(t)}(\check{c})=0$ holds, we have the minimum point is $\check{c}_i(t)=\frac{1}{t},~i=0,1,\cdots,t-1$.
It follows that
$$Reg_T(\hat{u}(t))\leq Reg_T(u_1(t)),$$
which yields that the online control policy $\hat{u}(t)$ serves as a better performer than $u_1(t)$.
\hfill $\Box$

\section{Proof of Lemma \ref{lemma02}}\label{app02}
\textbf{Proof.}
Consider the following quadratic optimization problem
\begin{align}
&minimize~~W_{T}=\sum_{t=0}^{T} z'(t)Qz(t)+v'(t)Rv(t), \nonumber \\
&subject~to~~~ z(t+1)=Az(t)+Bv(t).\label{lq01}
\end{align}
It follows from Theorem \ref{thm01} with $w(t)=0$ that the optimal index value of (\ref{lq01}) is
\begin{equation}
W_T^*=z(0)' P_T(0)z(0).
\end{equation}
Due to the time-invariance of the Riccati equation (\ref{ric01}), for any $0\leq t\leq T$, we have $P_T(t)=P_{T-t}(0)$.
For any $z(0)$ and $0\leq t_1<t_2\leq T$, it follows that
\begin{equation*}\begin{aligned}
z(0)'P_T(t_1)z(0)=W_{T-t_1}^*\geq W_{T-t_2}^*=z(0)'P_T(t_2)z(0),
\end{aligned}
\end{equation*}
which indicates that $P_T(t_1)\geq P_T(t_2)$.
Since $(A,B)$ is stabilizable, there exists a stabilizing control policy $v_s(t)=K_sz(t)$ such that $\lim_{t\rightarrow\infty} \|z_s(t)\|^2=0$ and
\begin{equation*}
\sum_{t=0}^{\infty}  \|z_s(t)\|^2\leq c_1 \|z(0)\|^2.
\end{equation*}
In this case, we have
\begin{equation*}\begin{aligned}
z(0)'P_T(t)z(0)=&~W_{T-t}^*\leq \sum_{i=0}^{T-t} z'(i)Qz(i)+v_s'(i)Rv_s(i) \\
\leq & \sum_{t=0}^{\infty} z'(t)Qz(t)+v_s'(t)Rv_s(t) \\
\leq & \sum_{t=0}^{\infty} z'(t)(Q+K_s'RK_s)z(t) \\
\leq & ~c\|z(0)\|^2,
\end{aligned}
\end{equation*}
where $c=\lambda_{\max}(Q+K_s'RK_s)c_1$. Thus, $P_T(t)$ is bounded.
Moreover, we have
$$\lim_{t\rightarrow-\infty}P_T(t)=\lim_{T\rightarrow\infty}P_T(0)=\hat{P},$$
where $\hat{P}$ satisfies the ARE (\ref{are01}).
Moreover, since $(A,B)$ is stabilizable and $(A,Q^{\frac{1}{2}})$ is observable, it follows from Theorem 1 in \cite{huang2008} that the ARE (\ref{are01}) has a unique positive definite solution $\hat{P}>0$. The proof is completed.
\hfill $\Box$

\section{Proof of Theorem \ref{thm04}}\label{app03}
\textbf{Proof.}
Since $D_T$ is independent with the control policy $u(t)$, we only need to consider the following quadratic optimization problem
\begin{align}
&minimize~~J_1(u(t))= \sum_{t=0}^{T} {\bf E} \big[ \hat{x}'(t)Q(t)\hat{x}(t)+u'(t)R(t)u(t) \big] \nonumber \\
&~~~~~~~~~~~~~~~~~~~~~~~~~~~~
+{\bf E} \big[ \hat{x}'(T+1)P_{T+1}\hat{x}(T+1)\big], \nonumber \\
&subject~to~~~ \hat{x}(t+1)=A\hat{x}(t)+Bu(t)+s(t). \label{op01}
\end{align}
For each time $t=0,1,\cdots,T$, define $\mathcal{G}(t)=G(t,u_a(t))$ with
\begin{align*}
G(t,u(t))= {\bf E} \big[ \hat{x}'(t)Q(t)\hat{x}(t)+u'(t)R(t)u(t)
+\mathcal{G}(t+1)\big].
\end{align*}
The terminal condition is given as
\begin{equation}
\mathcal{G}(T+1)={\bf E} \left[ \hat{x}'(T+1)P_{T+1}\hat{x}(T+1) \right].
\end{equation}
Next, we show that
\begin{align*}
\mathcal{G}(t)=&~{\bf E}\big[ \hat{x}'(t)P_T(t)\hat{x}(t)+2\hat{x}'(t)L_T(t)\mu_w\big] +\sum_{i=t}^T M(i),
\end{align*}
where
\begin{align}
M(t)=& -\mu_w'(P_T(t+1)+L_T(t+1))'B \Upsilon_T(t)^{-1} B'\nonumber \\
&\times (P_T(t+1)+L_T(t+1))\mu_w+2\mu_w' L_T(t+1)\mu_w \nonumber \\
&+2{\bf Tr} (A'P_{T}(t+1)B \Upsilon_T(t)^{-1} B'P_{T}(t+1)A \bar{Q}_v) \nonumber \\
&+{\bf Tr}(P_T(t+1)Q_w)+{\bf Tr}(P_T(t+1)\bar{Q}_v) \nonumber  \\
&-{\bf Tr}(A'P_T(t+1)A \bar{Q}_v). \label{mt01}
\end{align}

For each time $t=0,1,\cdots,T$, it follows that
\begin{equation*}
\begin{aligned}
&G(t,u(t)) \\
=&~{\bf E} \Big[ \hat{x}'(t)Q(t)\hat{x}(t)+u'(t)R(t)u(t)+(A\hat{x}(t)+Bu(t) \\
&+s(t))'P_{T}(t+1) (A\hat{x}(t)+Bu(t)+s(t)) \\
&+2(A\hat{x}(t)+Bu(t)+s(t))'L_{T}(t+1)\mu_w \Big]+\sum_{i=t+1}^T M(i)\\
=&~{\bf E} \Big[ (u(t)+\Upsilon_T(t)^{-1}(B'P_{t+1}A\hat{x}(t)+B'P_{T}(t+1)\mu_w \\
&+B'L_{T}(t+1)\mu_w))'\Upsilon_T(t)(u(t)+\Upsilon_T(t)^{-1} \\
&\times (B'P_{T}(t+1)A\hat{x}(t)+B'(P_{T}(t+1)+L_{T}(t+1))\mu_w)) \\
&-(B'P_{T}(t+1)A\hat{x}(t)+B'(P_{T}(t+1)+L_{T}(t+1))\mu_w)' \\
&\times \Upsilon_T(t)^{-1}(B'P_{T}(t+1)A\hat{x}(t)+B'P_{T}(t+1)\mu_w \\
&+B'L_{T}(t+1)\mu_w)
+\hat{x}'(t)(Q+A'P_{T}(t+1)A)\hat{x}(t) \\
&+ 2\hat{x}'(T)A'P_{T}(t+1)\mu_w +2\hat{x}'(t)A'L_T(t+1)\mu_w\\
&-2u'(T)B'P_{T}(t+1)AC^{-1}v(t) \Big]+2\mu_w'L_T(t+1)\mu_w \\
&+{\bf Tr}(P_{T}(t+1)Q_w) +{\bf Tr}(P_{T}(t+1)\bar{Q}_v))\\
&-{\bf Tr}(A'P_{T}(t+1)A\bar{Q}_v) +\sum_{i=t+1}^T M(i).
\end{aligned}
\end{equation*}
By utilizing the control policy $u_a(t)$ in (\ref{in02}), we have
\begin{equation*}
\mathcal{G}(t)= {\bf E}\left[ \hat{x}'(t)P_T(t)\hat{x}(t)+2\hat{x}'(t)L_T(t)\mu_w\right]+\sum_{i=t}^T M(i).
\end{equation*}
If we set $M_T=\sum_{t=0}^T M(t),$ it follows that
\begin{align*}
J_1(u_a(t))= &\sum_{t=0}^{T} {\bf E} \big[ \hat{x}'(t)Q(t)\hat{x}(t)+u_a'(t)R(t)u_a(t) \big] \nonumber \\
&+{\bf E} \big[ \hat{x}'(T+1)P_{T+1}\hat{x}(T+1)\big] \nonumber \\
=& ~\hat{x}'(0)P_T(t)\hat{x}(0)+2\hat{x}'(0)L_T(0)\mu_w
\\
&+M_T - D_T,  \nonumber
\end{align*}
where
\begin{equation*}
\begin{aligned}
&M_T-D_T \\
=&\sum_{t=0}^{T}\Big\{ -\mu_w'(P_T(t+1)+L_T(t+1))'B \Upsilon_T(t)^{-1} B'\\
&\times (P_T(t+1)+L_T(t+1))\mu_w+2\mu_w' L_T(t+1)\mu_w \\
&+2{\bf Tr}(A'P_{T}(t+1)B \Upsilon_T(t)^{-1} B'P_{T}(t+1)A \bar{Q}_v)\\
&+{\bf Tr}(P_T(t+1)Q_w)+{\bf Tr}(P_T(t+1)\bar{Q}_v)\\
&-{\bf Tr}(A'P_T(t+1)A \bar{Q}_v)
\hspace{-0.5mm}-\hspace{-0.5mm}
{\bf Tr}(Q(t)\bar{Q}_v)\hspace{-0.5mm}+\hspace{-0.5mm}{\bf Tr}(P_{T+1}\bar{Q}_v) \Big\}
\\
=&\sum_{t=0}^{T}\Big\{ -\mu_w'(P_T(t+1)+L_T(t+1))'B \Upsilon_T(t)^{-1} B'\\
&\times (P_T(t+1)+L_T(t+1))\mu_w+2\mu_w' L_T(t+1)\mu_w \\
&+{\bf Tr}(A'P_{T}(t+1)B \Upsilon_T(t)^{-1} B'P_{T}(t+1)A \bar{Q}_v)\\
&+{\bf Tr}(P_T(t+1)Q_w)\Big\}
-{\bf Tr}(P_T(0)\bar{Q}_v).
\end{aligned}
\end{equation*}
This proof is completed. \hfill $\Box$

\section{Proof of Theorem \ref{thm05}}\label{app04}
\textbf{Proof.}
For the optimization problem (\ref{op01}), define the following Lyapunov function
\begin{align*}
W(t)=&~{\bf E} \big[ \hat{x}'(t)P_T(t) \hat{x}(t)+2\hat{x}'(t)L_T(t)\mu_w \big]+\sum_{i=t}^T M(i),
\end{align*}
where $M(t)$ satisfies (\ref{mt01}).
It follows that
\begin{align}
&W(t)-W(t+1)\nonumber \\
=&~ {\bf E}  \left[  \hat{x}'(t)Q(t)\hat{x}(t)+u'(t)R(t)u(t) \right] \nonumber \\
&-{\bf E}  \left[ (u(t)-u_a(t))'\Upsilon_T(t)(u(t)-u_a(t))\right] \nonumber\\
&+2 {\bf E} \left[ u'(t)B'P_{T}(t+1)AC^{-1}v(t)\right] \nonumber\\
&+ 2{\bf Tr} \left(A'P_{T}(t+1)B \Upsilon_T(t)^{-1} B'P_{T}(t+1)A \bar{Q}_v\right)
\label{vt10}
\end{align}
where $u_a(t)$ is given in (\ref{in02}).
Summarizing (\ref{vt10}) from $t=0$ to $t=T$ yields that
\begin{align*}
&W(0)-W(T+1)\\
=&~\sum^T_{t=0}W(t)-W(t+1) \\
=&~ \hat{x}'(0)P_T(0)\hat{x}(0)+2\hat{x}'(0)L_T(0)\mu_w+\sum_{t=0}^T M(t) \\
& -{\bf E}[ \hat{x}'(T+1)P_T(T+1)\hat{x}(T+1)]\\
\end{align*}
\begin{align*}
=& \sum^T_{t=0} {\bf E}  \left[  \hat{x}'(t)Q(t)\hat{x}(t)+u'(t)R(t)u(t) \right] \\
&-\sum^T_{t=0}{\bf E}  \left[ (u(t)-u_a(t))'\Upsilon_T(t)(u(t)-u_a(t))\right]\\
&+ \sum^T_{t=0}2 {\bf E} \left[ u'(t)B'P_{T}(t+1)AC^{-1}v(t)\right] \\
&+ \sum^T_{t=0} 2{\bf Tr} \left(A'P_{T}(t+1)B \Upsilon_T(t)^{-1} B'P_{T}(t+1)A \bar{Q}_v\right),
\end{align*}
which implies that
\begin{align}
&J_1(u(t))\nonumber \\
=&~ \hat{x}'(0)P_T(0)\hat{x}(0)+2\hat{x}'(0)L_T(0)\mu_w+\sum_{t=0}^T M(t)  \nonumber \\
&+\sum^T_{t=0}{\bf E}  \left[ (u(t)-u_a(t))'\Upsilon_T(t)(u(t)-u_a(t))\right] \nonumber\\
&- \sum^T_{t=0} 2{\bf Tr} \left(A'P_{T}(t+1)B \Upsilon_T(t)^{-1} B'P_{T}(t+1)A \bar{Q}_v\right)
\nonumber \\
&-\sum^T_{t=0}2 {\bf E} \left[ u'(t)B'P_{T}(t+1)AC^{-1}v(t)\right] .
\end{align}
By utilizing the admissible control policy $\hat{u}_a(t)$ in (\ref{in01}),
we obtian

\begin{align*}
&J_T(\hat{u}_a(t))=J_1(\hat{u}_a(t))-D_T \\
=&~ \hat{x}'(0)P_T(t)\hat{x}(0)+2\hat{x}'(0)L_T(0)\mu_w+M_T-D_T  \nonumber \\
&+\sum^T_{t=0}{\bf E}  \left[ (\hat{u}_a(t)-u_a(t))'\Upsilon_T(t)(\hat{u}_a(t)-u_a(t))\right] \nonumber\\
&-\sum^T_{t=0}2 {\bf E} \left[ \hat{u}_a(t)B'P_{T}(t+1)AC^{-1}v(t)\right]  \nonumber \\
&- \sum^T_{t=0} 2{\bf Tr} \left(A'P_{T}(t+1)B \Upsilon_T(t)^{-1} B'P_{T}(t+1)A \bar{Q}_v\right) \\
=& ~ \hat{x}'(0)P_T(t)\hat{x}(0)+2\hat{x}'(0)L_T(0)\mu_w+H_T  \nonumber \\
&+{\bf Tr} \left({D}_T(0) \mu_w \mu_w'  \right)
+\sum^{T}_{t=1} \frac{1}{t} {\bf Tr}\left( {D}_T(t)C_w \right).
\end{align*}
Moreover, by Theorem \ref{thm04}, the regret satisfies
\begin{align}
\bar{R}eg_T(\hat{u}_a(t))=& ~J_T(\hat{u}_a(t))-J_T({u}_a(t))  \nonumber \\
=&\sum^{T}_{t=1} \frac{1}{t} {\bf Tr}\left( {D}_T(t)C_w \right).
\end{align}
Under hypotheses H1)-H2), we obtain that $\bar{R}eg_T(\hat{u}_a(t)) \leq O(\ln T)$. \hfill $\Box$






\end{document}